\documentclass[12pt]{article}
\usepackage{amssymb, amsmath}

\newtheorem{theorem}{Theorem}[section]

\newtheorem{definition}{Definition}[section]

\newcommand{\n}{\nonumber}

\newcommand{\si}{\sigma_R (|x|)}

\newcommand{\s}{\sigma}

\newcommand{\bb}{\begin{equation}}
\newcommand{\ee}{\end{equation}}
\newcommand{\bq}{\begin{eqnarray}}
\newcommand{\eq}{\end{eqnarray}}
\newcommand{\bqn}{\begin{eqnarray*}}
\newcommand{\eqn}{\end{eqnarray*}}

\begin{document}
\title{ On the nonexistence of stationary weak solutions to the compressible fluid equations}
\author{ Dongho Chae\\
 Department of Mathematics\\
  Sungkyunkwan
University\\
 Suwon 440-746, Korea\\
 e-mail: {\it chae@skku.edu }}
 \date{}
\maketitle
\begin{abstract}
In this paper we prove that under some integrability conditions for
the density and the velocity fields the only stationary weak
solutions to the compressible fluid equations on $\Bbb R^N$
correspond to the zero density. In the case of compressible
magnetohydrodynamics equations similar integrability conditions for
density, velocity and the magnetic fields lead to the zero density
and the zero magnetic field.
{\bf keywords:} compressible fluid equations,  compressible MHD
equations, nonexistence of stationary solutions
\end{abstract}
\section{Introduction}
 \setcounter{equation}{0}
 \subsection{The compressible Navier-Stokes(Euler) equations}
We are concerned  here on the compressible Navier-Stokes(Euler for
$\mu=\lambda =0$) equations on $\Bbb R^N$, $N\geq 1$.
  $$
(NS)\left\{\aligned &\partial_t \rho + \mbox{div}(\rho v) = 0, \\
& \partial_t (\rho v) + \mbox{div}(\rho v \otimes v) = -\nabla
p+\mu \Delta v +(\mu +\lambda ) \nabla \mathrm{div }\, v +f, \\
& p = p(\rho, S)\geq 0, \quad p=0 \,\, \mbox{only if}\,\,\rho=0.
\endaligned \right.
$$
The system (NS) describes  compressible gas flows with the adiabatic
exponent $\gamma$,  and  $\rho, v, S$, $p$ and $f$ denote the
density, velocity, specific entropy, pressure and the external force
respectively. Since the results below does not depend on the
particular form of the entropy equation, nor the specific form of
$p(\rho, S)$, we omit specifications of them.
 We treat  the viscous case $\mu
>0$(compressible Navier-Stokes equations) and the inviscid case
$\mu=\lambda=0$(compressible Euler equations) simultaneously.  For
surveys of the known mathematical theories of the equations we refer
to \cite{che, chr, lio, fei, nov}. Our aim here is to prove
nonexistence of nontrivial stationary weak solutions to the system
(NS) under suitable integrability conditions. In the inviscid case
and viscous case with $2\mu +\lambda =0$ our integrability condition
covers the finite energy condition, while for the viscous case with
$2\mu+\lambda \neq 0$  we need extra integrability for velocity
$v\in L^{\frac{N}{N-1}} (\Bbb R^N )$ besides the finite energy
condition. These could be regarded as Liouville theorem for the
stationary compressible fluid equations. The Liouville type of
theorems for the {\em nonstationary} incompressible Euler equations
are  recently studied by the author of this paper in \cite{cha1,
cha2}, where we need to impose extra condition for the sign of the
integral of pressure as well as the integrability conditions for the
velocity. In the case of compressible fluid equations, however, we
do not need such extra sign condition for pressure integral, since
the sign of pressure is automatically nonnegative pointwise. Similar
nonexistence results hold for the compressible magnetohydrodynamic
equations for $N\geq2$, which will be treated in the next section. A
stationary weak
 solutions of  (NS)  are defined as follows.
  \begin{definition} We say that  a triple $(v,\rho, S)\in [L^2_{loc} (\Bbb R^N)]^N\times
  L^\infty_{loc}(\Bbb R^N )\times  L^\infty_{loc} (\Bbb R^N )$ is a  stationary weak solution of $(NS)$
  if
  \bq\label{11}
 && \int_{\Bbb R^N} \rho v\cdot \nabla \psi \,dx=0 \qquad \forall \psi
 \in C_0 ^\infty (\Bbb R^N),\\
 \label{12}
 &&\int_{\Bbb R^N} \rho v\otimes v :\nabla \phi \,dx=-\int_{\Bbb R^N} p\, \,\mathrm{div }\, \phi
 \,dx-\mu\int_{\Bbb R^N} v\cdot \Delta \phi\, dx\n \\
 &&\qquad -(\mu +\lambda
 )\int_{\Bbb R^N} v\cdot \nabla \mathrm{div}\, \phi \, dx
 -\int_{\Bbb R^N} f\cdot \phi \, dx \quad \forall \phi
 \in [C_0 ^\infty (\Bbb R^N)]^N,\n\\
 \ \\
 \label{14}
 &&p = p(\rho, S)\geq 0, \quad p=0 \,\,\mbox{only if}\,\, \rho=0.
  \eq
  \end{definition}
  The following is our main nonexistence theorem for (NS).
 \begin{theorem} Let $N\geq 1$, and let the external force $f\in [L^1_{loc}(\Bbb R^N)]^N$
 satisfy $\mathrm{div}\, f=0$ in the sense
 of distribution. Suppose  $(\rho, v, S)$ is a
 stationary weak solution to $(NS)$ satisfying one of the following
 conditions depending on $\mu$ and $\lambda$.
 \begin{itemize}
 \item[(i)] In the inviscid case($\mu=\lambda =0$); there exists $w\in L^1_{loc} ([0, \infty))$,
 which is positive almost everywhere on $[0, \infty)$ such that
\bb\label{15} \int_{\Bbb R^N} \frac{(\rho |v|^2 +p)}{1+|x|}\, dx
 <\infty.
 \ee
 \item[(ii)] In the viscous case($\mu >0$);
 \begin{description}
 \item[(a)] if $2\mu+\lambda=0$,
 \bb\label{15b}
 \int_{\Bbb R^N} (\rho |v|^2 +p)\,dx
 <\infty.
\ee
 \item[(b)]
 if  $2\mu+\lambda\neq0$,
\bb\label{15c}
  \int_{\Bbb R^N} (\rho |v|^2
+|v|^{\frac{N}{N-1}}+p)\,dx
 <\infty.
 \ee
 \end{description}
 \end{itemize}
 Then, $\rho(x)=0$ for almost every $x\in \Bbb R^N$.
\end{theorem}
{\em Remark 1.1 } We note that in the special case of $ N\geq 3$,
$\mu>0, \mu+\lambda
>0$, $f=\rho \nabla \Phi$, where the connected component of $ \{
\Phi(x) >-c \}$ is unbounded, P.L. Lions showed nonexistence of
stationary solutions under appropriate integrability condition for
$\rho, v$(see Section 6.7 of \cite{lio}). Even when $\mu>0$,
$\mu+\lambda
>0$ and
$N\geq 3$, the above theorem  does not have mutual implication
relation with this.
\subsection{The compressible  MHD equations}
Here we are concerned  on  the compressible magnetohydrodynamic
equations on $\Bbb R^N$,
  $$
(MHD)\left\{\aligned &\partial_t \rho + \mbox{div}(\rho v) = 0, \\
& \partial_t (\rho v) + \mbox{div}(\rho v \otimes v -H\otimes H) =
-\nabla (p+\frac12 |H|^2)\\
&\qquad\qquad\qquad\qquad+\mu \Delta v +(\mu +\lambda ) \nabla \mathrm{div }\, v +f, \\
 &\partial_t H-\mathrm{curl}\, (v\times H)=0,\\
 & \mathrm{div}\, H=0,\\
& p = p(\rho, S)\geq 0, \quad p=0 \,\,\mbox{only if}\,\, \rho=0.\\
\endaligned \right.
$$
The system (MHD) describes compressible charged gas flows(plasma
gas) with the adiabatic exponent $\gamma$, and $\rho, v, H, S, p$
and $f$ denote the density, velocity, magnetic field, specific
entropy, pressure and the external force respectively. A {\em
stationary} weak
 solution of  (MHD)
 is defined as follows.
  \begin{definition} We say that a triple $(\rho, v, H, S)\in L^\infty_{loc}(\Bbb R^N )\times [L^2_{loc} (\Bbb R^N)]^N\times[L^2_{loc} (\Bbb R^N)]^N
  \times
    W^{1,2}_{loc}(\Bbb R^N )$ is a stationary weak solution of $(MHD)$
  if
  \bq\label{21}
 &&\int_{\Bbb R^N} \rho v\cdot \nabla \psi (x)\,dx=0 \qquad \forall \psi
 \in C_0 ^\infty (\Bbb R^N),\\
 &&\int_{\Bbb R^N} (\rho v\otimes v -H\otimes H ) :\nabla\phi (x)\,dx \n\\
 &&\qquad =-\int_{\Bbb R^N}( p+\frac12|H|^2)\,
  \mathrm{div }\, \phi (x)\,dx-\mu\int_{\Bbb R^N} v\cdot \Delta \phi\, dx\n \\
  \n\\
 &&\label{22}\qquad -(\mu +\lambda
 )\int_{\Bbb R^N} v\cdot \nabla \mathrm{div}\, \phi \, dx
 -\int_{\Bbb R^N} f\cdot \phi \, dx
 \quad \forall \phi
 \in [C_0 ^\infty (\Bbb R^N)]^N,\n\\
 \ \\
 &&\label{23}
 \int_{\Bbb R^N} (v\times H)\cdot \mathrm{curl}\, \varphi(x)\,dx=0\quad \forall \varphi
 \in [C_0 ^\infty (\Bbb R^N)]^N,\\
 &&\label{24}
  \int_{\Bbb R^N} H\cdot \nabla \eta (x)\, dx=0\qquad \forall
\eta
 \in C_0 ^\infty (\Bbb R^N)\\
 \label{26}
 &&p = p(\rho, S)\geq 0, \quad p=0 \,\,\mbox{only if}\,\, \rho=0.
  \eq
  \end{definition}
  Similarly to Theorem 1.1 we have the following theorem for (MHD).
 \begin{theorem} Let the external force $f\in [L^1_{loc}(\Bbb R^N)]^N$
 satisfy $\mathrm{div}\, f=0$ in the sense
 of distribution.  Suppose $(\rho, v,H, S)$ is a  stationary weak solution to
 (MHD) satisfying the following conditions depending on $\mu$ and $\lambda$.
  \begin{itemize}
  \item[(i)] In the inviscid case($\mu=\lambda=0$);
 \begin{description}
 \item[(i-a)] The case $N\geq 3$ : There exists $w\in L^1_{loc} ([0, \infty))$,
 which is non-increasing, positive almost everywhere on $[0,
 \infty)$ such that
\bb\label{27} \int_{\Bbb R^N} \frac{\rho |v|^2 +|H|^2+p}{1+|x|}\, dx
 <\infty
 \ee
 \item[(i-b)] The case $N=2$ :
 \bb\label{27aa}
\int_{\Bbb R^N} (\rho |v|^2 +|H|^2+p)dx
 <\infty.
 \ee
 \end{description}
\item[(ii)] In the viscous case($\mu >0$) for all $N\geq2$ ;
 \begin{description}
 \item[(ii-a)] if $2\mu+\lambda=0$,
\bb\label{27a} \int_{\Bbb R^N} (\rho |v|^2 +|H|^2+p)\, dx
 <\infty.
 \ee
 \item[(ii-b)] if $2\mu+\lambda\neq0$,
 \bb\label{27b}
\int_{\Bbb R^N} (\rho |v|^2 +|v|^{\frac{N}{N-1}}+|H|^2+p)\, dx
 <\infty.
 \ee
 \end{description}
 \end{itemize}
 Then, $\rho=0$ and $H=0$(and $v=0$ in the case (ii-b))  almost everywhere on $ \Bbb R^N$ if $N\geq
 3$.
 In the case $N=2$ we just conclude that $\rho=0$ almost
 everywhere on $\Bbb R^N$.
\end{theorem}
\noindent{\em Remark 2.1 } Contrary to the case of previous section,
our argument of the proof of the above theorem does not work for
$N=1$, and we do not yet know if similar nonexistence results hold
or not in those
cases.\\

\section{Proof of the Main Theorems}
\setcounter{equation}{0} In order to prove Theorem 1.1 we introduce
a class of weight functions as follows.
 \begin{definition}
We say  that a function $w(\cdot)\in C^3 ([0,\infty)) $ is
admissible if it satisfies the following conditions:
\begin{itemize}
\item[(i)]
\bb\label{w1} w(r), w'(r), w^{\prime\prime}(r)\geq 0\quad \mbox{and}
\quad w^{\prime\prime\prime}(r)\leq 0 \quad \forall r> 0.
  \ee
 \item[(ii)] There exists a constant $C$ such that
 \bb\label{w2}
  w^{\prime\prime} (r)+\frac{1}{r}w'(r) +\frac{1}{r^2}w(r)\leq
\frac{C}{1+r}\quad \forall r\geq 0.
 \ee
 \end{itemize}
 The class of all admissible function will be denoted by
 $\mathcal{W}$.
\end{definition}
As examples we find that $w_1 (r), w_2(r)\in \mathcal{W}$, where
$$ w_1 (r)=\log (\cosh r ), \quad w_2 (r)=\int^r_0 \arctan s\, ds. $$
\noindent{\bf Proof of Theorem 1.1 }\\
\noindent{\em \underline{(i) The case  $\mu=\lambda=0$ : }}
 Let us
consider a radial cut-off function $\sigma\in C_0 ^\infty(\Bbb R^N)$
such that
 \bb\label{18}
   \sigma(|x|)=\left\{ \aligned
                  &1 \quad\mbox{if $|x|<1$}\\
                     &0 \quad\mbox{if $|x|>2$},
                      \endaligned \right.
 \ee
and $0\leq \sigma  (x)\leq 1$ for $1<|x|<2$.  Let us choose a weight
function $w\in \mathcal{W}$. Then, for each $R
>0$, we define
 \bb\label{110}
\varphi_R (x)=w(|x|)\s \left(\frac{|x|}{R}\right)=w(|x|)\s_R
(|x|)\in C_0 ^\infty (\Bbb R^N).
 \ee
We   choose the vector test function
 $\phi$ in (\ref{12}) as
 \bb\label{111}
  \phi= \nabla \varphi_R (x).
 \ee
 Then, after routine computations, the equation (\ref{12}) becomes
  \bq\label{112}
 \lefteqn{0=\int_{\Bbb R^N} \rho (x) \left[W^{\prime\prime} (|x|)
\frac{(v\cdot x)^2}{|x|^2} + w^{\prime} (|x|)
\left(\frac{|v|^2}{|x|}
-\frac{(v\cdot x)^2}{|x|^3}\right) \right] \si  \,dx }\hspace{.5in}\n\\
&&\quad+\int_{\Bbb R^N} \rho(x) w' (|x|) \s'
\left(\frac{|x|}{R}\right)
\frac{(v\cdot x)^2}{R|x|^2} \,dx\n \\
&&\quad+ \int_{\Bbb R^N} \frac{1}{R}\rho(x)\left( \frac{|v|^2}{|x|}
-\frac{(v\cdot
x)^2}{|x|^3} \right) \s'\left(\frac{|x|}{R}\right)w(|x|) \,dx\n\\
  &&\quad+\int_{\Bbb R^N} \rho(x)\frac{(v\cdot x)^2}{
R^2|x|^2} \s^{\prime\prime} \left(\frac{|x|}{R}\right) w(|x|)
\,dx \n \\
 &&\quad+ \int_{\Bbb R^N}p(x)\left[ w^{\prime\prime}
 (|x|) +(N-1)\frac{w' (|x|)}{|x|}\right]\sigma_R (|x|)
  \, dx\n \\
&&\quad+ \frac{2}{R}\int_{\Bbb R^N}p(x) w' (|x|)
 \s' \left(\frac{|x|}{R}\right)\, dx\n \\
&&\quad+ \frac{N-1}{R}\int_{\Bbb R^N}p(x)\frac{1}{|x|}\s'
\left(\frac{|x|}{R}\right) w(|x|) \, dx\n \\
&&\quad+ \int_{\Bbb R^N}p(x)\frac{1}{R^2} \s^{\prime\prime}
\left(\frac{|x|}{R}\right)w(|x|) \,
dx \n \\
 &&:=I_1+\cdots +I_8
\eq From the condition  (\ref{w2}) we find that
 \bq\label{114}
\lefteqn{\int_{\Bbb R^N} (\rho(x)|v(x)|^2 +|p(x)|)
\left[w^{\prime\prime}
 (|x|) +\frac{1}{|x|}w' (|x|) +\frac{1}{|x|^2} w(|x|)\right]
 dx}\hspace{1.in}\n \\
 &&\leq C\int_{\Bbb R^N} \frac{\rho(x)|v(x)|^2 +|p(x)|}{1+|x|}\, dx
 < \infty.
 \eq
Since
\bqn
 &&\int_{\Bbb R^N} \rho(x)\left|\left[w^{\prime\prime} (|x|)
\frac{(v\cdot x)^2}{|x|^2} + w^{\prime} (|x|)
\left(\frac{|v|^2}{|x|} -\frac{(v\cdot x)^2}{|x|^3}\right)
\right]\right|dx\\
&&\qquad \leq 2
  \int_{\Bbb R^N} \rho(x)|v(x)|^2\left[ w^{\prime\prime}
(|x|)+\frac{w'(|x|)}{|x|} \right] dx <\infty, \eqn
  One can use the dominated convergence theorem to show that
  \bb\label{115}
  I_1 \to \int_{\Bbb R^N} \rho(x)\left[w^{\prime\prime} (|x|)
\frac{(v\cdot x)^2}{|x|^2} + w^{\prime} (|x|)
\left(\frac{|v|^2}{|x|} -\frac{(v\cdot x)^2}{|x|^3}\right) \right]
\,dx
 \ee
  as $R\to \infty$.
Similarly,
  \bb\label{116}
  I_5\to \int_{\Bbb R^N}p(x)\left[ w^{\prime\prime}
 (|x|) +(N-1)\frac{w' (|x|)}{|x|}\right]
 \, dx
\ee as $R\to \infty$.
 For $I_2$ we estimate
 \bq\label{117}
  |I_2 |&\leq& \int_{R<|x|<2R} \rho(x)|v(x)|^2\left|\s'
\left(\frac{|x|}{R}\right)\right|
  \frac{w'(|x|)}{|x|} \frac{|x|}{R}dx\n \\
  &\leq &2 \sup_{1<s<2} |\s'(s)|
\int_{R<|x|<2R}\rho(x) |v(x)|^2
  \frac{w'(|x|)}{|x|} dx\n \\
\to 0
  \eq
 as $R\to \infty$ by the dominated convergence theorem.
Similarly
  \bq\label{118}
   |I_3|&\leq &2 \int_{R<|x|<2R}  \frac{|x|}{R} \rho(x) |v(x)|^2
\left|\s'\left(\frac{|x|}{R}\right)\right|\frac{w(|x|)}{|x|^2} \,dx \n \\
&\leq &4\sup_{1<s<2} |\s'(s)|
 \int_{R<|x|<2R}\rho(x) |v(x)|^2
  \frac{w'(|x|)}{|x|}dx
\to 0,\n \\
  \eq
  and
  \bq\label{119}
  |I_4|&\leq&\int_{R<|x|<2R}\frac{|x|^2}{R^2} \rho(x)|v(x)|^2\left|\s^{\prime\prime}
  \left(\frac{|x|}{R}\right)\right|\frac{w(|x|)}{|x|^2} \,dx\n \\
  &\leq&4\sup_{1<s<2} |\s^{\prime\prime}(s)|
 \int_{R<|x|<2R}\rho(x)|v(x)|^2
  \frac{w(|x|)}{|x|^2}\,dx \to 0\n \\
  \eq
   as $R\to \infty$. The estimates for $I_6,I_7$ and $I_8$ are
   similar to the above, and we find
   \bq\label{120}
   |I_6|&\leq &2 \int_{R<|x|<2R}|p(x)| \frac{|x|}{R}\frac{w'
   (|x|)}{|x|}
 \left|\s' \left(\frac{|x|}{R}\right)\right|  \, dx\n \\
 &\leq& 4\sup_{1<s<2} |\s' (s)|
\int_{R<|x|<2R}|p(x)|\frac{w'
   (|x|)}{|x|}dx \to 0,\n \\
  \eq
\bq\label{121}
   |I_7|&\leq & (N-1)\int_{R<|x|<2R}|p(x)|\frac{|x|}{R}\left|\s'
\left(\frac{|x|}{R}\right)\right| \frac{w(|x|)}{|x|^2} \, dx\n \\
 &\leq& 2\sup_{1<s<2} |\s' (s)|
\int_{R<|x|<2R}|p(x)|\frac{w(|x|)}{|x|^2}dx\to 0,\n \\
  \eq
  and
  \bq\label{122}
 |I_8|&\leq&\int_{\Bbb R^N}|p(x)|\frac{|x|^2}{R^2}
\left|\s^{\prime\prime}\left(\frac{|x|}{R}\right)\right|
\frac{w(|x|)}{|x|^2} \, dx \n\\
  &\leq& 4\sup_{1<s<2}
|\s^{\prime\prime} (s)|
 \int_{R<|x|<2R}|p(x)|\frac{w
   (|x|)}{|x|^2}dx \to 0\n \\
  \eq
  as $R\to \infty$ respectively. Thus passing $R\to \infty$ in (\ref{112}),
  we finally obtain
\bq\label{123}
 &&\int_{\Bbb R^N} \rho(x)\left[w^{\prime\prime} (|x|)
\frac{(v\cdot x)^2}{|x|^2} + w^{\prime} (|x|)
\left(\frac{|v|^2}{|x|} -\frac{(v\cdot x)^2}{|x|^3}\right) \right]
 \,dx \n \\
&&\qquad+\int_{\Bbb R^N}p(x)\left[ w^{\prime\prime}
 (|x|) +(N-1)\frac{w' (|x|)}{|x|}\right] \, dx=0.\n \\
 \eq
Since
$$w^{\prime\prime} (|x|)
\frac{(v\cdot x)^2}{|x|^2} + \frac{w^\prime (|x|)}{|x|}\left(|v|^2
-\frac{(v\cdot x)^2}{|x|^2}\right)\geq 0,
$$
 and
 $$
 w^{\prime\prime}
 (|x|) +(N-1)\frac{ w^\prime (|x|)}{|x|}>0,
 $$
 in (\ref{123}), we need to have
 \bqn
 p(x)=p(\rho(x), S(x))=0 \qquad\mbox{almost everywhere on $\Bbb
 R^N$},
 \eqn
and
 therefore $\rho(x)=0$ for almost every $x\in \Bbb R^N$.\\
 \ \\
 \noindent{\underline{\em (ii) The case of $\mu >0$ and either $2\mu +\lambda =0$ or  $2\mu +\lambda \neq0$ : }}   In this case
 we choose the vector test function
 $$\phi =\nabla (|x|^2\sigma_R (x))
 $$
 in
 (\ref{12}), where $\sigma_R$ is defined above. Then, each of the
 procedure of (i) can be repeated word by word with specific choice of function $w(r)\equiv
 1$ on $[0, \infty)$. We just need to show
 \bb\label{125}
 \mu\int_{\Bbb R^N} v\cdot \Delta \phi dx +(\mu +\lambda)\int_{\Bbb
 R^N} v\cdot \nabla \mathrm{div}\, v dx=o(1)
 \ee
 as $R\to \infty$.\\
 If $2\mu+\lambda =0$, then
\bqn
  J&:=&\mu\int_{\Bbb R^N} v\cdot \Delta \nabla (|x|^2 \s_R )\,
dx +(\mu +\lambda)\int_{\Bbb
 R^N} v\cdot \nabla [\mathrm{div}\, \nabla (|x|^2 \s_R)]\, dx\\
 &=&2(\mu +\lambda)\int_{\Bbb R^N} v\cdot\nabla \Delta (|x|^2
 \s \left(\frac{|x|}{R}\right)\, dx=0,
 \eqn
and (\ref{125}) holds true. \\
If $2\mu+\lambda \neq0$, then we
compute and estimate
 \bqn |J|&=&2|\mu +\lambda| \left|\int_{\Bbb R^N} v\cdot\nabla \Delta (|x|^2
 \s \left(\frac{|x|}{R}\right)\, dx\right|\\
&\leq& |2\mu +\lambda| \left|\int_{\Bbb R^N}
(N+5)\left[\frac{(v\cdot x)}{R|x|}\s ' \left(\frac{|x|}{R}\right)+
\frac{(v\cdot x)}{R^2}\s^{\prime\prime}
\left(\frac{|x|}{R}\right)\right]\, dx\right|\n \\
&&\qquad+|2\mu
+\lambda|\left|\int_{\Bbb R^N} \frac{|x|(v\cdot
x)}{R^3}\s^{\prime\prime\prime}
\left(\frac{|x|}{R}\right)\,dx \right|\n \\
&&\leq  \frac{C}{R} \int_{R\leq |x|\leq 2R} |v(x)| \, dx \leq C
\left(\int_{R\leq  |x|\leq 2R} |v(x)|^{\frac{N}{N-1}}
dx\right)^{\frac{N-1}{N}} \to 0
 \eqn
   as $R\to \infty$, since $v\in L^{\frac{N}{N-1}} (\Bbb R^N)$ by the hypothesis in the viscous case.
    Thus (\ref{125})  holds true.
 $\square$\\
 \ \\
\noindent{\bf Proof of Theorem 1.2 }  Since the proof is similar to
that of Theorem 1.1, we will be brief here.  Let $w\in \mathcal{W}$.
In the inviscid case we choose the vector test function
 $\phi$ in (\ref{22}) as previously, namely
 \bb\label{29}
  \phi= \nabla \varphi_R (x),
 \ee
 where
 $$\varphi_R(x)=w(|x|)\s_R (|x|)=w(|x|)
 \s\left(\frac{|x|}{R}\right)
 $$
 and $\sigma(\cdot)$ is the cutoff function given by (\ref{18}).
 Then, the equation (\ref{22}) become
\bq\label{210}
 &&\int_{\Bbb R^N} \rho(x) \left[w^{\prime\prime} (|x|) \frac{(v\cdot
x)^2}{|x|^2} + \frac{w^{\prime} (|x|)}{|x|}  \left(|v|^2
-\frac{(v\cdot x)^2}{|x|^2}\right)\right]
 \s_R (x)\,dx \n \\
&&-\int_{\Bbb R^N} \left[w^{\prime\prime} (|x|) \frac{(H\cdot
x)^2}{|x|^2} + \frac{w^{\prime} (|x|) }{|x|}  \left(|H|^2
-\frac{(H\cdot x)^2}{|x|^2}\right)\right]
 \s_R (x)\,dx \n \\
&&\qquad+\int_{\Bbb R^N}(p(x)+\frac12 |H|^2)\left[ w^{\prime\prime}
 (|x|) +(N-1)\frac{w^{\prime} (|x|)}{|x| } \right] \s_R
 (x)\,dx\n \\
 &&=o(1)
 \eq
 as $R\to \infty$. Passing $R\to \infty$ in (\ref{210}), and rearranging
 the remaining terms, we have
\bq\label{211}
  &&\int_{\Bbb R^N} \rho(x) \left[w^{\prime\prime}  (|x|) \frac{(v\cdot
x)^2}{|x|^2} + \frac{w^{\prime} (|x|)}{|x|} \left(|v|^2
-\frac{(v\cdot x)^2}{|x|^2}\right)\right]
\,dx \n \\
&&\quad+\int_{\Bbb R^N} \left[ \frac{w^{\prime}
(|x|)}{|x|}-w^{\prime\prime} (|x|)\right]
\frac{(H\cdot x)^2}{|x|^2}\, dx\n \\
 &&\quad+\frac{N-3}{2}\int_{\Bbb R^N} \frac{|H|^2 w^{\prime}
(|x|)}{|x|}\,dx+\frac12\int_{\Bbb R^N}  |H|^2 w^{\prime\prime}(|x|)\, dx \n \\
&&\qquad+\int_{\Bbb R^N}p(x)\left[w^{\prime\prime}
 (|x|) +(N-1)\frac{w^{\prime}
(|x|)}{|x|} \right]
 \,dx=0\n \\
 \eq
 Since $w^{\prime\prime\prime}(r)\leq 0$, we find that $w^{\prime\prime}(r)$
 is a nonnegative, non-increasing function on $ [0, \infty)$. Hence,
 $$ \frac{w^{\prime} (|x|)}{|x|}=\frac{w^{\prime\prime}(|x|)+\int_0 ^{|x|} w^{\prime\prime} (s)\, ds}{|x|} \geq
 \frac{\int_0 ^{|x|} w^{\prime\prime} (s)\, ds}{|x|}\geq w^{\prime\prime}
(|x|),
$$
and
$$ \frac{w^{\prime}
(|x|)}{|x|}-w^{\prime\prime} (|x|)\geq 0
$$
 for almost every $x\in \Bbb R^N$.  Therefore for $N\geq 3$,
all of the integrals in (\ref{211}) are nonnegative, and hence
$$p=p(\rho, S)=0,\qquad H=0, $$
and $\rho=0, H=0$ almost everywhere on $\Bbb R^N$. For the case
$N=2$ we just set $w(r)= r^2 $ on $[0, \infty)$. Then, (\ref{211})
is reduced to
$$
\int_{\Bbb R^N} [\rho(x) |v(x)|^2 + 2p(x)]dx=0,
$$
which implies $p(\rho, S)=0$, and hence $\rho=0$ almost everywhere
on $\Bbb R^N$.
 The proof for the
viscous case is the same as that of Theorem 1.1, and we omit it
here.
$\square$\\
$$\mbox{ \bf Acknowledgements} $$
The author would like to thank to Seung-Yeal Ha  and G-Q. Chen for
useful discussion. This work was supported partially by  KRF
Grant(MOEHRD, Basic Research Promotion Fund)

\end{document}